\documentclass[12pt]{article}
\usepackage{amsmath, amsfonts, amsthm, amssymb, epsfig}

\newtheorem{theorem}{Theorem}
\newtheorem{lemma}{Lemma}

\begin{document}
\title{{\bf Reciprocal Properties of Pythagorean triangles}}
\author{Konstantine Zelator\\
Mathematics, Statistics, and Computer Science\\
212 Ben Franklin Hall\\
Bloomsburg University of Pennsylvania\\
400 East 2nd Street\\
Bloomsburg, PA  17815\\
USA\\
and\\
P.O. Box 4280\\
Pittsburgh, PA  15203\\
kzelator@.bloomu.edu\\
e-mails: konstantine\underline{\ }zelator@yahoo.com}

\maketitle

\section{Introduction} {\it Crux Mathematicorum with Mathematical Mayhem}, is
a problem solving journal published by the Canadian Mathematical Society for
undergraduate university students and teachers, as well as high school
teachers and students.  It publishes material ranging from K-10 to K-14 grade
levels.  The problems are divided into four categories: skoliad problems,
mayhem problems, mathematical olympiad problems, and general problems. In the
April, 2009 issue of the journal (Vol. 35, No; see reference \cite{1}), the
following problem was published (Mayhem problem M390), quote: ``A Pythagorean
triangle is a right-handed angled triangle with all three sides of integer
length.  Let $a$ and $b$ be the legs of a Pythagoream triangle and $h$ be
altitude to the hypotenuse.  Determine all triangles for which,

$$
\dfrac{1}{a} + \dfrac{1}{b} + \dfrac{1}{h} = 1 {\rm ''},\ {\rm end\ of\
  quote}.
$$

A solution was published in the February 2010 issue of the journal (see
\cite{2}).  It was shown that the only Pythagorean triangle with the above
property is the triangle with leg lengths $3$ and $4$; and hypotenuse length
$5$.  Our solution to the same problem appears in the proof of part (iii) of
Theorem 2, found in Section 5 of this article.

The above problem provides the motivation behind this work. Throughout, we
will use the standard notation $(u,w)$ to denote the greatest common divisor
of two integers $u$ and $v$.

The following definition creates conceptual basis for this work.

\vspace{.15in}

\noindent {\bf Definition 1:}  {\it Let $(a,b,c)$ be a Pythagorean triple with
  $a$ and $b$ being the leg lengths; $c$ the hypotenuse length.  Also let $h$
  be the altitude to the hypotenuse.  Let $v,k,l$ be positive integers with
  $k$ and $l$ being relatively prime; $(k,l) =1$.  Then, the Pythagorean
  triple is to have the reciprocal property $R(v,k,l)$ if the positive
  integers $a,b,h,v,k,l$ satisfy the condition or equation,

$$
\dfrac{1}{a} + \dfrac{1}{b} + \dfrac{v}{h} = \dfrac{k}{l}
$$}

\vspace{.15in}

According to Definition 1 and the solution to the above problem, the
Pythagorean triple that has the reciprocal property $(R(1,1,1)$, is the triple
$(3,4,5)$.

Let us outline the organization of this article.

In Section 2, we state the well known parametric formulas that generate the
entire set of Pythagorean triples.  We provide two references, \cite{3} and
\cite{4}.  Section 3 contains three results from number theory.  The first
one, Lemma 1, is very well known and it is generally known as Euclid's Lemma.
We use Lemma 1 to prove Lemma 2 which is in turn used in the proof for Theorem
4 (Section 5).  The third result in Section 3 is Theorem 1 which is very well
known.  It gives the general parametric solution to a linear diophantine
equation in two variables; which, in turn, is used in establishing Theorem 5.
Section 4 contains the key equation of this article.  The results of this work
are expressed in Theorems 2, 3, 4, and 5, all found in Section 5.

\section{Pythagorean triples}  

\noindent \parbox{5.0in}{{\bf Definition 2} {\it A triple $(a,b,c)$ of three positive integers $a,b,c$ is said to be Pythagorean with hypotenuse length $c$ if
  $a^2+b^2 = c^2$.

A positive integer triple $(a,b,c)$ is Pythagorean if, and only if, $a =
d(m^2-n^2),\ b = d(2mn),\ c = d(m^2+n^2)$; ($a$ and $b$ may be switched).
where $m,n,d$ are positive integers such that $m > n$ (and so $m \geq 2$)\ $(m,n) = 1$ (i.e., $m$ and $n$ are relative prime), and $m+n \equiv 1({\rm
  mod}\,2)$ (i.e., one of $m,n$ is even; the other odd)  If $d=1$, then triple
$(a,b,c)$ is said to be primitive.}} \hfill (1)

\section{Three results from number theory} 

For a proof of Lemma 1, see reference \cite{3}.

\begin{lemma}  (Euclid's lemma)  Suppose that $\alpha, \beta, \gamma$ are
  positive integers such that $\alpha$ divides the product $\beta \gamma$;
  with $\alpha$ and $\beta$ being relatively prime.  Then $\alpha$ is a
  divisor of $\gamma$.
\end{lemma}

\begin{lemma} Let $\dot{q}$ and $r>q$ be two relatively prime positive
  integers with different parities; one of them being even, the other
  odd. Also, assume that $q$ is not divisible by $3$.  That is, $q \not\equiv
  0({\rm mod}\,3)$ and that $r\not\equiv q({\rm mod}\,3)$. Then,

$$
\left( 3r^2 + q^2 + 2rq,\ \ 2rq (r^2-q^2)\right) = 1
$$
\end{lemma}

\begin{proof}  To establish the result, we will show that no prime divisor of
  the product $2rq(r-q)(r+q)$ divides the integer $3r^2 + q^2 + 2rq$.  First
  observe that by virtue of the fact that one of $r,q$ is odd, while the other
  is even.  The integer $3r^2+q^2 + 2rq$ is clearly odd.

Now, let $p$ be an odd prime divisor of the product $2rq(r-q)(r+q)$.  The
prime $p$ must divide one of the four factors $r,q,\ r-q,$ or $r+q$.  If $p$
divides $r$ and $p$ also divides $3r^2+q^2 + 2rq$, then it follows that $p$
must divide $q^2$; and since $p$ is a prime, this implies that $p$ is a
divisor of $q$.  So $p$ divides both $r$ and $q$ violating the condition
$(r,q)=1$, i.e., the hypothesis that $r$ and $q$ are relatively prime.  Now,
if $p$ divides $q$ and $p$ divides $3r^2 + q^2 + 2rq$, a straight forward
calculation implies that $p$ divides $3r^2$.  So, $p$ is a common divisor of
both $q$ and $3r^2$.  But $q$ is not divisible by $3$ so $p \neq 3.$ Since $p$
divides $3r^2$ and $(p,3) = 1$, Lemma 1 implies that $p$ must divide $r^2$,
and hence, $p$ divides $r$.  Once again, this is contrary to the hypothesis
$(r,q)=1$.  Next, suppose that $p$ divides $r-q$ or that $p$ divides $r+q$.
In other words, $r \equiv \pm q({\rm mod}\,p)$.  Consequently,

$$
\begin{array}{rcl}
3r^2 + q^2 + 2rq & \equiv & 4q^2 \pm 2q^2({\rm mod}\, p);\\
\\
3r^2 + q^2 + 2rq & \equiv & 6q^2 \ {\rm or}\ 2q^2({\rm mod}\,p)
\end{array}
$$

\noindent So, if $p$ also divides $3r^2 + q^2 + 2rq$, then, if the minus sign
holds in the above congruences, $p$ must also divide $2q^2$.  But $p$ is odd,
so this implies that $p$ divides $q^2$, and so $p$  divides $q$. But then $r
\equiv -q({\rm mod}\,p)$ implies that $p$ divides $r$ as well; contrary to the
condition $(r,q) =1$.  

Finally, if the plus sign holds in the above congruences, it follows that $p$
must divide $6q^2$.  If $p$ divides $q^2$, then we obtain the same
contradiction as before (in the minus sign case) which means (the only
possibility left) that $p$ must divide $6$; and thus (since $p$ is odd),
$p=3$.  But, then since the plus sign holds, we get $r \equiv q({\rm mod}\,3)$
contrary to the hypothesis that $r\not\equiv q({\rm mod}\,3)$.  The proof is
complete.  \end{proof}

The theorem below is very well known.  For a reference, see \cite{3}.

\begin{theorem} Let $A$ and $B$ be integers, not both zero, and $D$ their
  greatest common divisor, $D=(A,B)$, and $C$ an integer.  Consider the
  two variable linear diophantine equation $Ax+By = C$.

\begin{enumerate}
\item[(i)]  If $D$ is not a divisor of $C$, then the above equation has not
  integer solution.

\item[(ii)]  If $D$ is a divisor of $C$, then the entire integer solution set
  can be described  by the parametric formulas, $x=x_0 + \frac{B}{D} t,\ y =
  y_0 - \frac{A}{D} t$ where $t$ can be any integer and $\{x_0,y_0\}$ is a
  particular integer solution.
\end{enumerate}
\end{theorem}

\section{The key equation}

Let $(a,b,c)$ be a Pythagorean triple and $v,k,l$ given positive integers with
$(k,l) = 1$.    If the triple $(a,b,c)$ has the reciprocal property $R(v,k,l)$,
we then have 

\setcounter{equation}{1}
\begin{equation}
\dfrac{1}{a} + \dfrac{1}{b} + \dfrac{v}{h} = \dfrac{k}{l} \label{E2}
\end{equation}

\noindent where $h$ is the altitude to the hypotenuse.  We have Area $A =
\frac{1}{2} ab = \frac{1}{2} hc$, which gives

\begin{equation}
h = \dfrac{ab}{c} \label{E3}
\end{equation}

From (\ref{E2}) and (\ref{E3}) we obtain

$$
\dfrac{1}{a} + \dfrac{1}{b} + \dfrac{vc}{ab} = \dfrac{k}{l}
$$

\noindent which further gives

\begin{equation}
l(b+a+vc) = k \,ab. \label{E4}
\end{equation}

Combining (\ref{E4}) with the parametrical formulas in (1) yields

\begin{equation}
l \left[ (v+1)m^2+(v-1)n^2 + 2mn\right] = d\cdot k\cdot (2mn)(m^2-n^2)
\label{E5}
\end{equation}

\noindent which is the key equation.

\section{ Theorems 2, 3, 4, 5 and their proofs}

\begin{theorem}
\begin{enumerate}
\item[(i)] Let $(a,b,c)$ be a Pythagorean triple described by the parametric
  formulas in (1).  If $(a,b,c)$ has the reciprocal property $R(v,k,1)$, then
  it is necessary that $m$ is a divisor of $v -1$ and $n$ a divisor of $v+1$.

\item[(ii)]  There exists no Pythagorean triple that has the reciprocal
  property $R(2,k,1)$.

\item[(iii)]  (Our solution to problem M390) The only Pythagorean triple which
  has reciprocal property $R(1,1,1)$ is the triple $(3,4,5)$.

\item[(iv)]  If $k \geq 2$, there exists no Pythagorean triple which has the
  reciprocal property $R(1,k,l)$.

\item[(v)] Let $v,k,l$ be positive integers with $(k.l)=1$, $l$ odd and $v$
  even.  Then, there exists no Pythagorean triple which has the reciprocal
  property $R(v,k,l)$.
\end{enumerate}
\end{theorem}

\begin{proof}
\begin{enumerate} 
\item[(i)] Since $l=1$, the key equation (\ref{E5}) implies 

$$
(v+1)m^2 + (v-1)n^2 + 2mn = d\cdot k \cdot (2mn)(m^2-n^2)
$$

\noindent which gives 

$$(v-1)n^2 = m \left[ 2dkn (m^2-n^2) - 2n - m (v+1)\right].
$$

Since $(m,n) = 1 = (m,n^2)$, the last equation and Lemma 1 imply that $m$ must
be a divisor of $v-1$.  Likewise, from the above equation we get 

$$(v+1)m^2 = n\left[ 2dkm (m^2-n^2) - 2m - n(v-1)\right].
$$

A similar argumentm as in the previous case, establishes that $n$ must be a
divisor of $v+1$.

\item[(ii)] If $v=2$ and $l = 1$, then it follows from part (i) that $m$ must
  be a divisor of $v-1 = 2-1 = 1$ which is impossible since $m \geq 2$ (see
  (1)).

\item[(iii)] For $v = k = l =1$, equation (\ref{E5}) takes the form

\begin{equation}
\begin{array}{rcl}
2m^2 + 2mn & = & d(2mn)(m^2-n^2);\\
\\
2m(m+n) & = & d(2mn)(m-n)(m+n);\\
\\
1 & = & dn(m-n).
\end{array}\label{E6}
\end{equation}

Since $m,n,$ and $m-n$ are positive integers (in view of $m > n \geq 1$).  The
last equation (\ref{E6}) implies $d = n = m-n=1$ which gives $d=1,\ n=1,\
m=2$.  So, by (1) we obtain $(a,b,c) = (3,4,5)$

\item[(iv)]  For $v=1$, equation (\ref{E5}) takes the form

\begin{equation}
\begin{array}{rcl}
l(2m^2 + 2mn) & = & d \cdot k(2mn)(m^2-n^2);\\
\\
2ml(m+n) & = & d\cdot k (2mn)(m-n)(m+n);\\
\\
l & = & dk(m-n)
\end{array} \label{E7}
\end{equation}

Since $k$ is relatively prime to $l$ and $k \geq 2$, clearly, equation is
impossible or contradictory.

\item[(v)]  Since $v$ is even, both integers $v-1$ and $v+1$ are odd.  Since
  $m$ and $n$ have different parities (i.e., one is odd and the other even),
  it follows that $(v+1)m^2 + (v-1)n^2 + 2mn$ is an odd integer; and since $l$
  is odd.  The left-hand side of equation (\ref{E5}) is then odd, while the
  right-hand side is even, creating a contradiction.  So, no Pythagorean
  triple can have a reciprocal property $R(v,k,l)$ with $v$ being even and $l$
  being odd.
\end{enumerate}
 \end{proof}

\begin{theorem} Let $k$ be a positive integer and $v$ a positive integer such
  that both $v-1$ and $v+1$ are primes (i.e., twin primes).  Then, there
  exists no Pythagorean triangle or triple that has the reciprocal property
  $R(v,k,1)$.
\end{theorem}

\begin{proof}
If, to the contrary such a triangle or triple exists, then equation (\ref{E5})
with $l=1$ gives

\begin{equation}
(v+1)m^2 + (v-1)n^2 +2mn = d\cdot k\cdot (2mn)(m^2-n^2) \label{E8}
\end{equation}

We put $v-1=p$, a prime; and so $v+1=p+2$, also a prime.  By part (i), we know
that the positive integer $m$ must be a divisor of $v-1=p$; and since $p$ is a
prime and $m\geq 2$, it follows that $m=p$.  Likewise, from part (i), we know
that $n$ must divide $v+1 = p+2$.  But, $p+2$ is a prime.  So, either $n=1$ or
$n=p+2$ which cannot be the case since $m > n$, and so $p >n$.  Hence, the only
remaining possbility is $n=1$.

We have $m=p,\ n=1,\ v-1=p$, and $v+1=p+2$.  Accordingly, by equation
(\ref{E8}) we get 

$$
\begin{array}{rcl}
(p+2)p^2 + p + 2p & = & d\cdot 2\cdot p (p-1)(p+1)k;\\
\\
(p+2)p + 3 & = & 2dk(p^2-1).
\end{array}
$$

After some algebra 

\begin{equation}
(2dk-1)p^2 - 2p - (3+2dk) = 0 \label{E9}
\end{equation}

Equation (\ref{E9}) demonstrates that the prime $p$ is one of the two roots or
zeros of the quadratic trinomial $t(x) = (2dk-1) x^2-2x - (3+2dk)$.

Since this trinomial has integer coefficients and an integer root, the other
root must be a rational number (since the sum and the product of the roots are
both rational). This, in turn, necessitates that the discriminant must be a
perfect or integer square:

$$4+4(2dk-1)(3+2dk) = T^2,
$$

\noindent for some non-negative integer $T$. Obviously, $T$ must be even, $T=2t$ for some $t \in {\mathbb Z},\ t \geq 0$.

We obtain

\begin{equation}
\begin{array}{rcl}
1+(2dk-1)(3+2dk) & = & t^2;\\
\\
4dk(dk+1) - 2 & = & t^2.
\end{array} \label{E10}
\end{equation}

However, $t^2 \equiv 1$ or $0({\rm mod}\,4)$, depending on whether $t$ is odd
or even.  Hence, the left-hand side of (\ref{E10}) is congruent to $2$ modulo
$4$.  While the right-hand side is congruent to $1$ or $0$ modulo $4$, a
contradiction.  The proof is complete. \end{proof}

\begin{theorem}Let $(a,b,c)$ be a Pythagorean triangle or triple with the two
  (out of three) generator positive integers $m$ and $n$ in (1), satisfying
  the additional conditions $n \not\equiv 0({\rm mod}\,3)$ and $m \not\equiv
  n({\rm mod}\, 3)$.

Under this additional hypothesis the following holds true.  The triple
$(a,b,c)$ has the reciprocal property $R(2,1,l)$ if, and only if, the positive
integers $d$ (third generator in the formulas in (1)) and $l$ are of the
form,

$$
\begin{array}{rcll} d & = & t \cdot (2mn)(m^2-n^2)(3m^2+n^2 + 2mn) & {\rm
    and}\\
\\
l & = & t(2mn)(m^2-n^2) \end{array}
$$

\noindent where $t$ is a positive integer.\end{theorem}

\begin{proof}
For $v=2$ and $k=1$, the key equation (\ref{E5}) takes the form

\begin{equation}
l(3m^2 + n^2 + 2mn) = d(2mn)(m^2-n^2). \label{E11}
\end{equation}

A straight forward calculation shows that if $l$ and $d$ have the required
form, then equation (\ref{E11}) is satisfied.  Now, the converse.  Assume that
(\ref{E11}) is satisfied.  Since $n \not\equiv 0({\rm mod}\,3),\ m \not\equiv
n({\rm mod}\, 3)$, and also by (1), $m > n$ and one of $m,n$ is even, the other
odd.  We infer from Lemma 2 that the positive integer  $2mn(m^2-n^2)$ is
relatively prime to the positive integer $3m^2+n^2 +2mn$.  Since it divides
the product on the left-hand side of (\ref{E11}), we deduce that
$2mn(m^2-n^2)$ must be a divisor of $l$.

Hence, $l = t(2mn)(m^2-n^2)$, for some positive integer $t$.  Substituting for
$l = t(2mn)(m^2-n^2)$, in (\ref{E11}), produces $d=t (2mn)(m^2-n^2) (3m^2+n^2
+ 2mn)$. \end{proof}

\begin{theorem}
The Pythagorean triple $(3,4,5)$ has the reciprocal property $R(v,k.l)$
precisely when the three positive integers $v,k,$ and $l$ belong to one of six
groups.

$$\begin{array}{llll}
{\rm Group\ 1:} & v = 1+12t, & k=1+5t,  & l = 1 \\
{\rm Group\ 2:} & v = 1+6t,  & k=2+5t,  & l = 2\\
{\rm Group\ 3:} & v = 1+4t,  & k=3+5t,  & l = 3\\
{\rm Group\ 4:} & v = 1+3t,  & k=4+5t,  & l = 4\\
{\rm Group\ 5:} & v = 1+2t,  & k=6+5t,  & l = 6\\
{\rm Group\ 6:} & v = 1+t,   & k=12 +5t,& l=12
\end{array}
$$

\noindent where $t$ can be any non-negative integer in all six groups.
\end{theorem}

\begin{proof}  Going back to (1), we see that for the Pythagorean triple
  $(3,4,5)$ we have $d=1,\ m = 2$, and $n=1$.  Therefore, the key equation
  (\ref{E5}) yields,

\begin{equation} l(7+5v) = 12 k.  \label{E12}
\end{equation}

By virtue of $(l,k) = 1$, equation (\ref{E12}) and Lemma 1 imply that $l$ must
be a divisor of $12$.  Thus,

$$ l = 1,2,3,4,6, \ {\rm or}\ 12.
$$

\noindent For each of these six values of $l$, a linear diophantine equation
the the variables $v$ and $k$ is obtained.

We tabulate these below by including in each case an easily obtained (by
inspection) particular solution$\{v_0,k_0\}$

\vspace{.15in}
$
\begin{array}{|l|c|l|} \hline
& {\rm Linear\ diophantine\ equation} & {\rm Particular\ solution}\\ \hline
l = 1 & -5v+12 k=7 & v_0=1,\ k_0=1 \\ \hline
l = 2 & -5v+6k = 7 & v_0=1,\ k_0=2 \\ \hline
l = 3 & -5v+4k = 7 & v_0=1,\ k_0=3 \\ \hline
l = 4 & -5v+3k = 7 & v_0=1,\ k_0=4 \\ \hline
l = 5 & -5v+2k = 7 & v_0=1,\ k_0=6 \\ \hline
l = 6 & -5v+k = 7 & v_1=1,\ k_0=12 \\ \hline
\end{array}
$

\vspace{.15in}

Note that in each in the above table, the greatest common divisor of $D$ of
the coefficients is always $1$.  Applying Theorem 1 we obtain the following
parametric solutions.

\vspace{.15in}

$\begin{array}{ll}
{\rm For}\ l = 1, & v=1+12t,\ k=1+5t\\
{\rm For}\ l = 2, & v=1+6t,\ k=2+5t\\
{\rm For}\ l = 3, & v=1+4t,\ k=3+5t\\
{\rm For}\ l = 4, & v=1+3t,\ k=4+5t\\
{\rm For}\ l = 6, & v=1+2t,\ k=6+5t\\
{\rm For}\ l = 12, & v=1+t,\ k=12+5t
\end{array}
$

According to Theorem 1, $t$ can be any integer in the general solution.
However, in our case $v$ and $k$ are positive integers.  A cursory examination
shows that if $t$ is a negative integer in each of the above six cases, at
least one of $v$ and $k$ will be zero or negative.  Thus, the parameter $t$
cannot take any negative values.  On the other hand, in all six groupings, for
each non-negative value of $t$, the resulting values of $v$ and $k$ are both positive.
\end{proof}

\noindent Remark: If one cannot find, by inspection, a particular solution
  $\{x_0,y_0\}$ to a linear diophantine equation (Theorem 1).  There is the
  Euclidean algorithm that guarantees that such a solution can be
  found.

\end{document}